\def\AA         {{\bf A}}
\def\ZZ         {{\bf Z}}

\def\CC         {{\bf C}}

\def\PP         {{\bf P}}

\def\ee         {{\rm e}}
\def\no         {\hspace{2pt}:\hspace{-2pt}}
\def\on         {\hspace{-2pt}:\hspace{2pt}}

\def\vac        {|0\rangle}

\def\Fock       {{\rm Fock}}

\def\dim        {{\rm dim}}

\documentstyle[twoside,12pt]{article}
\newtheorem{prop}{Proposition}[section]
\newtheorem{dfn}[prop]{Definition}
\newtheorem{theo}[prop]{Theorem}

\newtheorem{rem}[prop]{Remark}

\title{Introduction to the vertex algebra approach to mirror symmetry}

\author{
Lev A. Borisov
\\
\small Department of Mathematics,  Columbia University \\
\small 2990 Broadway, Mailcode 4432, New York, NY 10027, USA \\
\small e-mail: lborisov@math.columbia.edu}

\begin{document}

\date{}

\maketitle

\begin{abstract}
{
The goal of this paper is to make the vertex operator algebra approach
to mirror symmetry accessible to algebraic geometers. Compared to better-known
approaches using moduli spaces of stable maps and special Lagrangian
fibrations, this approach follows more closely the original line of thinking 
that lead to the discovery of mirror symmetry by physicists. The ultimate 
goal of the vertex algebra approach is to give precise mathematical 
definitions of N=(2,2) superconformal field theories called A and B models
associated to any Calabi-Yau variety and then show that thus constructed 
theories are related by the mirror involution for all known examples of 
mirror symmetric varieties.
}
\end{abstract}

%\section{plan}
%I will eventually delete this.
%
%1. Introduction.
%
%2. General review of superconformal field theories. Superconformal
%   field theories on the sphere and vertex algebras.
%
%3. Definition and basic properties of vertex algebras. OPEs. Definitions
%   of conformal structure, N=2 structure, topological structure. 
%   BRST operators.
%
%4. Examples of vertex algebras. One free boson. One free fermion.
%   Free bosons and fermions based on a pair of dual spaces. Vertex
%   algebras based on the lattice. Vertex algebra of a fan.
%
%5. MSV construction. (Quasi)-loop-coheherent sheaves mentioned.
%
%6. Summary of the results of my big paper. A short remark about
%   the paper with Libgober. 
%
%7. Open questions.

\section{Introduction}
This paper should serve as an introduction to the vertex algebra
approach to mirror symmetry developed in \cite{borvert}. It is
thus understandable that our emphasis and selection of topics
reflects the author's bias. The reader should keep in mind that 
other approaches to mirror symmetry exist and have independent
mathematical interest. In particular, the stable maps approach 
allowed to state and prove mathematically the predictions for the 
(virtual) numbers of rational curves on a quintic threefold and
other similar examples, see \cite{kontsevich,givental,yauetal}.

It is widely stated in physics literature that given a Calabi-Yau
manifold $X$ together with an element of its complexified K\"ahler cone
one can construct two N=(2,2) superconformal field theories 
called $A$ and $B$ models, see for example \cite{witten}.
On the other hand, all actual calculations and definitions of these 
theories involve Feynman type integrals
over infinite-dimensional spaces of all maps from Riemann surfaces
to $X$. While physicists have developed a good intuitive understanding
of the formal properties of these integrals, they are mathematically
ill-defined.

The precise axiomatic definition of $N=(2,2)$ superconformal field theory 
that would include the $A$ and  $B$ models above is still not available.
Roughly speaking, this theory is a modular functor, see for example
\cite{Segal}, but the number of labels is, perhaps, infinite. 
In particular, there must exist a Hilbert space $H$ such that every 
Riemann surface whose oriented boundary consists of $k$ incoming 
and $l$ outcoming circles produces an operator from $H^{\otimes k}$ to
$H^{\otimes l}$, perhaps defined only up to a scalar multiple.
Superconformal field theory is a highly complicated object. Even when 
the Riemann surface is a sphere, the structure of superconformal field 
theory is rather non-trivial. A typical way to construct such a theory 
is by building it from the representation theory of vertex algebras 
that satisfy certain conditions, see for example  \cite{huang}. In fact, 
mirror symmetry originated from the work of Gepner \cite{gepner} who 
used (finite quotients of) tensor products of the so-called minimal 
models which are certain irreducible representations 
of the N=2 superconformal algebra. He was able to match the dimensions of 
chiral rings (see \cite{LVW}) of the resulting theories with the dimensions 
of the cohomology spaces of the Calabi-Yau hypersurfaces in projective 
spaces, in particular quintic threefolds.

Vertex algebra approach to mirror symmetry attempts to define rigorously
superconformal field theories associated to Calabi-Yau manifolds 
and then prove that the corresponding theories for mirror manifolds 
are related to each other. At this stage only the vertex algebra of 
the theory has been recovered and much work is still to be done.
This review contains no new results, and no proofs are presented.
It is intended as an introduction to vertex algebras
for algebraic geometers, and its ultimate goal is to enable an
interested reader to understand the paper \cite{borvert}. In particular,
only vertex algebras that appear in the context of hypersurfaces in 
toric varieties are discussed.

Section 2 contains basic definitions and properties of vertex algebras,
and follows closely the book of Kac \cite{Kac}. Section 3 provides
the reader with a simplest non-trivial example of vertex algebra called one
free boson. It is generalized to several
bosons and several fermions in Section 4. Section 5 is devoted to 
the very important paper of Malikov, Schechtman and Vaintrob \cite{MSV}
who construct chiral de Rham complex of an arbitrary smooth variety.
We are mostly interested in the case of Calabi-Yau varieties.
Section 6 explains main results of \cite{borvert} and the last section
summarizes briefly the problems that are still to be addressed in
the vertex algebra approach.

This article is based in part on the talks given at Columbia, Northwestern, 
MIT and Rutgers. The author wishes to thank these institutions for their 
hospitality. The author also thanks  Ezra Getzler, Yi-Zhi Huang and 
Peter Landweber for useful references and stimulating conversations.

\section{Definition and basic properties of vertex algebras}
The goal of this section is to state definitions of vertex algebras 
and to introduce important notions of normal ordered products and
operator product expansions (OPE). Our treatment follows closely the 
book of Kac \cite{Kac}. We also define N=2 superconformal structures
and describe BRST cohomology construction necessary to understand
\cite{borvert}.

%definition of vertex algebra
\begin{dfn}
{\rm
(\cite{Kac})
A vertex algebra is the set of data that consists of a super vector space 
$V$ (over $\CC$), a {\it state-field correspondence} $Y$ and 
a {\it vacuum vector}\ $\vac$. The fact that $V$ is a superspace
simply means that $V=V_0\oplus V_1$. Elements of $V_0$ are called bosonic or
even and elements from $V_1$ are called fermionic or odd. Vacuum vector
$\vac$ is a bosonic element of $V$. The  most important structure is 
the state-field correspondence $Y$ which is a parity preserving linear map from
$V$ to ${\rm End} V[[z,z^{-1}]]$
$$a\,\line(0,1){5}\hspace{-3.6pt}\to Y(a,z)=\sum_{n\in Z}a_{(n)}z^{-n-1}$$
such that for every two elements $a$ and $b$ the elements 
$a_{(n)}b$ are zero for all sufficiently big $n$.
To form a vertex algebra the data $(V,Y,\vac)$ must satisfy the 
following axioms.
\\
\noindent$\bullet${\bf translation covariance:}
$\{T,Y(a,z)\}_-=\partial_z Y(a,z)$ where $\{,\}_-$ denotes
the usual commutator and $T$ is defined by $T(a)=a_{(-2)}\vac$;
\\
\noindent$\bullet${\bf vacuum:} $Y(\vac,z)={\bf 1}_V,~
Y(a,z)\vac_{z=0}=a$;
\\
\noindent$\bullet${\bf locality:} $(z-w)^N\{Y(a,z),Y(b,w)\}_{\mp}=0$
for all sufficiently big $N$, where $\mp$ is $+$ if and only if
both $a$ and $b$ are fermionic. The equality is understood as an identity
of formal power series in $z$, $z^{-1}$, $w$ and $w^{-1}$. It is often
expressed by saying that $Y(a,z)$ and $Y(b,z)$ are mutually local.
}
\end{dfn}

%Normal ordered products and OPEs
Let $a$ and $b$ be two elements of the vertex algebra $V$. We 
denote the corresponding {\it fields}\  $Y(a,z)$ and $Y(b,z)$ by $a(z)$
and $b(z)$ respectively. Locality axiom of the vertex algebra allows one 
to express the supercommutators of the {\it modes}\ $a_{(m)}$ and $b_{(n)}$ 
in a concise way in terms of {\it operator product expansions} (OPEs).
Namely, define {\it normal ordered product} $\no a(z)b(w)\on\in {\rm End}V
[[z,z^{-1},w,w^{-1}]]$ by the formula
$$
\no a(z) b(w) \on = \sum_{m\in \ZZ_{<0},n\in \ZZ} a_{(m)}b_{(n)}
z^{-m-1} w^{-n-1} \pm \sum_{m\in \ZZ_{\geq 0},n\in \ZZ} b_{(n)}a_{(m)}
z^{-m-1}w^{-n-1}
$$
where $\pm$ is $-$ if and only if both $a$ and $b$ are fermionic.
Then it is not hard to show (see \cite{Kac} for details) that
locality axiom implies 
$$
a(z)b(w)=\sum_{j=0}^{N-1}\frac {c^j(w)}{(z-w)^{j+1}}
+ \no a(z)b(w)\on
$$
where $c^j(w)$ are some elements of ${\rm End}V[[w,w^{-1}]]$ 
and $(z-w)^{-j-1}$ is Laurent expanded in the region $|z|>|w|$. 
Moreover, there holds a remarkable Borcherds OPE formula, which states that 
$c^j(w) = Y(a_{(j)}b,w)$
so $c^j$ are also fields in the vertex algebra.
All information about supercommutators of the modes of $a$ and $b$
is conveniently encoded in the $\sum$ part of this OPE and the parities of 
$a$ and $b$. The $\sum$ part is called {\it singular part} of the OPE.

%Conformal structure 
\begin{def}\label{confstruct}
{\rm A vertex algebra with a conformal structure is a vertex algebra
$(V,Y,\vac)$ with a choice of an even element $v$ such that the corresponding
field $Y(v,z)=:L(z)$ satisfies the operator product expansion
$$
L(z)L(w)= \frac{c/2}{(z-w)^4} + \frac{2L(w)}{(z-w)^2}
+\frac{\partial_wL(w)}{z-w} + \no L(z)L(w) \on
$$
where $c$ is a constant called {\it central charge}. In addition, one
assumes that $L_{(0)}$ coincides with the operator $T$ in the definition
of vertex algebra. We also assume that $L_{(-1)}$ is diagonalizable 
on $V$, all of its eigenvalues are real numbers, and 
$$
\{L_{(-1)},Y(a,z)\}_-=z\partial_zY(a,z)+Y(L_{(-1)}a,z)
$$
for all $a$. 
}
\end{def}

\begin{rem}
{\rm
The same vertex algebra can have many different conformal structures.
Element $v$ that defines a conformal structure is called Virasoro element.
}
\end{rem}

Once a conformal structure is fixed, it is customary to shift the 
index in the definition of $a_{(n)}$ as follows. If $a$ has eigenvalue
$\alpha$ with respect to $L_{(-1)}$ then we introduce the notation
$$
Y(a,z)=\sum_{n\in \ZZ}a_{(n)}z^{-n-1} =: \sum_{n\in \ZZ-\alpha}
a[n]z^{-n-\alpha}.
$$
The number $\alpha$ is called the {\it conformal weight} of $a$.
In particular, we observe that OPE of $L(z)L(w)$ implies that 
conformal weight of the Virasoro element is two, and we introduce
$L(z)=\sum_{n\in \ZZ}L[n]z^{-n-2}$. In these notations the endomorphisms
$L[m]$ satisfy the commutator relations 
$$
\{L[m],L[n]\}_- = (m-n)L[m+n]+\frac c{12}(m^3-m)\delta_{m+n}^0
$$
of the Virasoro algebra with central charge $c$.

%N=2 structure
When one studies Calabi-Yau manifolds, one 
obtains vertex algebras which have not only conformal structure, but
what is called N=2 superconformal structure. It consists of the 
choice of conformal structure plus an even field $J$ and two
odd fields $G^+$ and $G^-$ which satisfy the following OPE.
$$
L(z)L(w)= \frac{c/2}{(z-w)^4} + \frac{2L(w)}{(z-w)^2}
+\frac{\partial_wL(w)}{z-w} + \no L(z)L(w) \on,
$$
$$
L(z)J(w)=\frac{J(w)}{(z-w)^2} + \frac {\partial_w J(w)}{z-w} + 
\no L(z)J(w)\on,
$$
$$
L(z)G^{\pm}(w)=\frac{(3/2)G^{\pm}(w)}{(z-w)^2} + \frac{\partial_w G^{\pm}
(w)}{z-w}+ \no L(z)G^{\pm}(w)\on,
$$
$$
J(z)J(w)=\frac{c/3}{(z-w)^2}+\no J(z) J(w)\on,
$$
$$
J(z)G^{\pm}(w)=\pm\frac{G^{\pm}(w)}{z-w}+ \no J(z)G^{\pm}(w)\on,
$$
$$
G^{\pm}(z)G^{\mp}(w)=\frac{2c/3}{(z-w)^3}\pm\frac
{2J(w)}{(z-w)^2}
+\frac{2L(w)\pm\partial_wJ(w)}{z-w}+ \no G^{\pm}(z)G^{\mp}(w)\on ,
$$
$$
G^{\pm}(z)G^{\pm}(w)=\no G^{\pm}(z)G^{\pm}(w)\on.
$$
It is common to introduce N=2 charge $\hat c=c/3$, where
$c$ is the central charge of the usual Virasoro algebra.
%
%topological algebras
Also one often changes the notations slightly by introducing Virasoro
field $L_{top}=L(z)+(1/2)\partial_zJ(z)$ of conformal charge $0$.
Then $G^\pm$, $J$ and $L_{top}$ form {\it topological structure of dimension}\ 
$\hat c=d$. 

%mirror involution
Notice that if one switches $G^+$ and $G^-$ and changes the sign
of $J$, then one obtains another N=2 structure. This involution
is called {\it mirror involution} and it is expected to switch
$A$ and $B$ models constructed from mirror symmetric Calabi-Yau manifolds. 
(More precisely, it is supposed to act this way on the holomorphic part
and act by identity on the antiholomorphic part of N=(2,2) superconformal
field theory, but this paper only deals with the holomorphic part. 
In fact, the absence of precise understanding of how to put together 
holomorphic and anti-holomorrhpic parts of the theory is a big 
obstacle in the vertex algebra approach to mirror symmetry.)

%BRST operators
Given a vertex algebra $(V,Y,\vac)$ one can construct other algebras
by the BRST cohomology construction. If $a\in V$ is such that $a_{(0)}^2=0$,
then one considers cohomology of $V$ with respect to $a_{(0)}$, called
{\it BRST cohomology}. Operator $a_{(0)}$ and field $Y(a,z)$
are called BRST operator and BRST field respectively.
One can show, see for example \cite{borvert}, that 
BRST cohomology of $V$ with respect to $a_{(0)}$ has a
natural structure of vertex algebra. Moreover, if $a_{(0)}$ supercommutes
with the fields of N=2 structure, then this structure descends to
BRST cohomology.

\section{First example of vertex algebra: one free boson}
The simplest non-trivial example of the vertex algebra is called {\it one 
free (chiral) bosonic field}. The goal of this section is to describe
explicitly the data $(V,Y,\vac)$. Moreover, it turns out that this 
vertex algebra could be provided with a conformal structure, which we
also describe. Most of the calculations are skipped, and the reader 
is referred to \cite{Kac}.

Consider an abstract unital associative algebra generated by elements
$d[n], n\in \ZZ$ with commutator relations
$$\{d[m],d[n]\}_- = m\delta_{m+n}^0$$
In other words, $d[n]$ commutes with everything except $d[-n]$, and
commutators of $d[n]$ and $d[-n]$ are proportional to the identity.

There is a standard representation of this algebra called a {\it Fock space}.
Namely, consider a vacuum vector $\vac$ such that 
$$d[\ZZ_{\geq 0}]\vac = 0$$
and try to see what space could be built from it.
We will call operators $d[\ZZ_{>0}]$ {\it annihilators}\ and 
operators $d[\ZZ_{<0}]$ {\it creators}.  Operator $d[0]$ commutes with
all other operators and will equal zero on this Fock space. Its 
relevance will be seen later in Section 6 when we talk about 
vertex algebras defined by a lattice. 

Notice that all creators commute with each other. 
If we apply all creators to the vacuum vector, assuming that the results
are linearly independent, we get the space 
$$V=\oplus_{n_1,n_2,...}
\CC \prod_{k>0}d[-k]^{n_k} \vac$$
where all $n_k$ are non-negative integers, and only finitely many of them
are nonzero. Creators obviously act on this space. The action of
annihilators could be defined by means of the commutator rules and the
fact that annihilators vanish on the vacuum vector. For example,
$$d[3]\,d[-1]^2d[-3]^2d[-5]\vac 
=d[-1]d[-3]^2d[-5]d[3]\vac 
+ d[-1]\{d[3],d[-3]^2\}_-d[-5]\vac$$
$$
=6d[-1]d[-3]d[-5]\vac .$$
Thus $V$ is a representation of the algebra of $d$ and
it is called the Fock space of one free bosonic field. One can think of it
as the space of polynomials in infinitely many variables
$d[-1],d[-2],...$ with creators acting by multiplication and annihilators
acting by differentiation. 

To describe the structure of vertex algebra on this Fock space $V$, we
will need additional notations. We introduce the field $d(z)\in{\rm
End}V[{z,z^{-1}}]$ by the formula
$$d(z)=\sum_{n\in \ZZ} d[n] z^{-n-1}.$$
Notice that $d[z]\vac = \sum_{k\geq 0}d[-k-1]z^k$, and when you plug in
$z=0$ you get $d[-1]\vac$. Eventually $d[z]$ will be a field that
corresponds to $d[-1]\vac$. To construct other fields we use the notion
of normal ordering introduced in Section 2. If we try to make sense of 
$d(z)^2=d(z)d(z)$ as an element of ${\rm End} V[[z,z^{-1}]]$, we run
into infinities. On the other hand, one can plug  $w=z$ into 
$\no d(z)d(w)\on$ and the resulting field $\no d(z)d(z)\on$ makes sense
as an element of ${\rm End} V[[z,z^{-1}]]$. 
In terms of the modes, one can define
$$: d[m]d[n] :~=~\left\{ 
\begin{array}{ll}
d[m]d[n]&{\rm if}~ n\geq 0\\
d[n]d[m]&{\rm if}~ n<0\\
\end{array}
\right.
$$
and then write
$$:d(z)d(z):~=~\sum_{m,n\in \ZZ} :d[m]d[n]: z^{-m-n-2}.$$
We notice that this field applied to vacuum lies in $V[[z]]$ and 
$$:d(z)d(z):\vac_{z=0}= d[-1]^2\vac .$$
Similarly, one defines fields 
$$\no\prod_{k\geq 0}
\left(\frac{\partial^k d}{\partial z^k} \right)^{n_k}
\on$$
by pushing all annihilators to the right and all creators to the left.
Then the claim is that these fields form the fields from the definition
of the vertex algebra that correspond to the states
$$
\prod_{k\geq 0}
(k!)^{n_k}d[-k-1]^{n_k}\vac .
$$
One can also show that these fields are mutually local. 
This allows us to define the state-field correspondence $Y$ which 
satisfies vacuum and locality axioms.
One can also show that translation axiom is satisfied. Moreover,
the operator $T$ could be written in terms of $d[n]$ as
$$
T=\frac12\sum_{k\in \ZZ} d[k]d[-1-k].
$$
As a result, we have constructed our first example of the vertex algebra.

We will now describe how to equip this vertex algebra $(V,Y,\vac)$
with a conformal structure. Look at the field
$$L(z)=\frac 12 :d(z)d(z):$$
and introduce $L[n]$ by $L(z)=\sum_{n\in \ZZ} L[n]z^{n-2}$.
Explicitly these operators could be written as  $L[n]=(1/2)\sum_{k\in
Z}:d[k]d[n-k]: $. Then one can check that the following commutator 
relations hold
$$\{L[m],L[n]\}_- = (m-n)L[m+n] + {1\over 12}(m^3-m)\delta_{m+n}^0.$$
Therefore, $L[m]$ form Virasoro algebra with central charge one.

One observes that $T=L[-1]$ is the translation covariance operator 
$T$.  Also, $L[0]$ is diagonalizable, because 
$$L[0]d[-1]^{n_1}d[-2]^{n_2}... \vac=
(\sum_i {in_i})~d[-1]^{n_1}d[-2]^{n_2}... \vac.$$

We can conveniently rewrite the commutators of $d[n]$ in terms of 
the OPEs. After an easy calculation, we get 
$$
d(z)d(w)=\frac 1{(z-w)^2}  + \no d(z)d(w)\on.
$$
In general, it is straightforward to calculate OPEs
of products of free bosons and their derivatives using Wick's theorem,
see \cite{Kac}. The key point of Wick's theorem is that the commutator
of products of two sets of linear operators such that the pairwise
commutators are in the center can be explicitly written in terms of 
these pairwise commutators.

\section{Further examples of vertex algebras: several free fermions 
and bosons}
The construction of the previous section can be generalized in several 
directions. First of all, instead of considering one free boson, one
may consider several of them. This simply means that one considers
the tensor product of a number of copies of the Fock space of one
free boson with the structure of vertex algebra induced in an obvious
way. More generally, for every finite-dimensional vector space $W$
of dimension $r$ over $\CC$ equipped with a non-degenerate symmetric 
bilinear form  $\langle~,~\rangle$, one constructs a vertex algebra.
One considers a unital associative
algebra generated by $w[n]$, $n\in \ZZ$, $w\in W$ with the commutator
relations 
$$\{w_1[m],w_2[n]\}_- = m\langle w_1,w_2 \rangle \delta_{m+n}^0.$$
Then, analogously to the one-dimensional example, one defines a Fock space
generated from $\vac$ by applying negative modes of $w[n]$. This
Fock space carries a natural structure of vertex algebra, which is 
isomorphic to the tensor product of $\dim W$ copies of one free boson.

We will be particularly interested in the case of the space $W$ which 
is a direct sum of a space $W_1$ and its dual, and thus has a natural 
non-degenerate symmetric bilinear product denoted by $\cdot$. 
In terms of the OPEs, the algebra  is generated by the fields $a(z)$ 
and $b(z)$ where $a\in W_1$, $b\in W_1^*$ and the OPEs are 
$$
a_1(z)a_2(w)=\no a_1(z)a_2(w)\on,~~ b_1(z)b_2(w)=\no b_1(z)b_2(w)\on,
$$
$$
a(z)b(w) = \frac {a\cdot b}{(z-w)^2}+\no a(z)b(w) \on.
$$

So far we have not introduced any fermionic elements. This is easily
accomplished by changing commutators $\{~,~\}_-$ to anticommutators
$\{~,~\}_+$ in the above formulas. While the construction could be
described for a single free fermion, we will restrict our attention
to $2r$ free fermions constructed from $W_1\oplus W_1^*$ where
$W_1$ is a vector space of dimension $r$.  One starts with 
a unital associative algebra generated by $\varphi[n]$ and $\psi[n]$
with $\varphi\in W_1^*,~\psi\in W_1,~n\in \ZZ+\frac 12$ with the 
anticommutator relations 
$$
\{\varphi_1[m] ,\varphi_2 [n] \}_+ = 0;~\{\psi_1[m] ,\psi_2 [n] \}_+ = 0;
$$
$$
\{\varphi[m] ,\psi [n] \}_+ = (\varphi\cdot\psi) \delta_{m+n}^0.
$$
The Fock space is constructed by applying pairwise anticommuting creators
$\varphi[(\ZZ+\frac12)_{<0}],~\psi[(\ZZ+\frac12)_{<0}]$ to the vacuum
vector $\vac$ which is annihilated by the rest of the modes. 
The vertex algebra structure is provided by the products of various
derivatives of the fields
$$\varphi(z) = \sum_{n\in \ZZ+\frac 12} \varphi[n] z^{-n-\frac12},~~
\psi(z) = \sum_{n\in \ZZ+\frac 12} \psi[n] z^{-n-\frac12},$$
and the parity is defined by the total number of $\varphi$ and $\psi$.
Operator product expansions of the fields $\varphi(z)$ and $\psi(z)$ 
are 
$$
\varphi_1(z)\varphi_2(w)=\no \varphi_1(z)\varphi_2(w)\on,~~ 
\psi_1(z)\psi_2(w)=\no \psi_1(z)\psi_2(w)\on,
$$
$$
\varphi(z)\psi(w) = \frac {\varphi\cdot \psi}{z-w}+
\no \varphi(z)\psi(w) \on.
$$
One can provide this vertex algebra with the structure of conformal
vertex algebra by introducing a field 
$$L(z)=\frac 12 \no \partial \psi^i(z)\varphi_i(z)-\psi^i(z)\partial 
\varphi_i(z)\on$$ 
where $\{\psi^i\}$ and $\{\varphi_i\}$ are dual bases of $W_1$ and $W_1^*$.
We implicitly sum over all $i$ via standard physical convention.
A fermionic analog of Wick's theorem allows us to calculate that 
central charge of this conformal structure is $\dim W_1$.

Finally, we put together fermions and bosons. The resulting algebra
is a crucial component of the constructions of \cite{MSV} and \cite{borvert}.
Again, let $W_1$ and $W_1^*$
be two dual spaces. We consider the vertex algebra which is a product
of the free fermionic and free bosonic algebras constructed above.
It is generated by fields $\varphi(z)$, $\psi(z)$, $a(z)$ and $b(z)$.
It is equipped with a conformal structure of central charge
$c=3\dim W_1$. Moreover, one can extend this structure to an N=2 structure
$$
G^+(z):=\varphi_i (z) a^i(z),~
G^-(z):= \psi^i(z) b_i(z),~ 
J(z):= \no \varphi_i(z)\psi^i(z)\on,
$$
$$
L(z):=\frac 12 \no a^i(z)b_i(z)\on + 
\frac 12 \no \partial \psi^i(z)\varphi_i(z)-\psi^i(z)\partial \varphi_i(z)
\on
$$
where again $\{a^i\}$ and $\{b_i\}$ are dual bases of $W_1$ and $W_1^*$.
We remark that the resulting N=2 fields are independent from the choice
of these bases. The N=2 central charge is $\hat c = c/3 = \dim W_1$.

\section{Chiral de Rham complex}
In a breakthrough paper \cite{MSV} Malikov, Schechtman and Vaintrob
have introduced a sheaf of vertex algebras which they call chiral
de Rham complex for every complex manifold. Roughly speaking, the idea
is to associate to every manifold $X$ a sheaf which locally over a
neighborhood of a point $x\in X$ looks like a vertex algebra with 
$2\dim X$ bosons and $2\dim X$ fermions associated to the vector space
$W=T_X(x)\oplus T^*_x(X)$.

There are some important details that we need to address. First of all,
one needs to use a slightly different version of the definition 
of the vertex algebra of free bosons. Namely, instead of the commutator
relations
$$\{a[m],b[n]\}_- = m(a\cdot b)\delta_{m+n}^0$$
they use the relations
$$\{a[m],b[n]\}_- = (a\cdot b)\delta_{m+n}^0.$$
The vacuum is now annihilated by $a[\ZZ_{\geq 0}]$ and $b[\ZZ_{>0}]$,
but $b[0]$ are considered to be creators. The fields $a(z)$ and $b(z)$ 
are now defined by 
$$b(z):=\sum_{n\in \ZZ}b[n]z^{-n}$$
and the basic OPE is 
$$
a(z)b(w)=\frac {a\cdot b}{z-w} +\no a(z)b(w)\on.
$$ 
Roughly speaking, one uses $\int b(w) dw$ instead of $b(w)$. 
The fields of the N=2 algebra are modified accordingly.
$$
G^+_{\rm MSV}(z):= \varphi_i (z) a^i(z),~
G^-_{\rm MSV}(z):= \psi^i(z) \partial_z b_i(z),~ 
J_{\rm MSV}(z):= \no \varphi_i(z)\psi^i(z)\on,
$$
$$
L_{\rm MSV}(z):=\frac 12 \no a^i(z)\partial_z b_i(z)\on + 
\frac 12 \no \partial \psi^i(z)\varphi_i(z)-\psi^i(z)\partial \varphi_i(z)
\on.
$$
The crucial observation of \cite{MSV} is that the group of 
automorphisms of the ring of local coordinates of $X$ at $x$ embeds
into the group of vertex algebra automorphisms of the above vertex algebra.
This allows one to glue together the above spaces and construct a sheaf of 
vertex algebras over the variety $X$. Unfortunately, the fields of the
$N=2$ algebra are not preserved under general automorphisms. However, if
$X$ is a Calabi-Yau variety, then the existence of the holomorphic volume
form allows one to restrict the attention to the volume-preserving local
changes of coordinates, and these changes do preserve the fields of $N=2$ 
algebra.

In \cite{MSV} the resulting sheaf is called {\it chiral de Rham complex},
because the usual de Rham complex is naturally embedded in it.
We will denote the chiral de Rham complex by  ${\cal MSV}(X)$. 
An important remark here is that it is not a quasi-coherent sheaf. 
The multiplication map ${\cal O}(X)\times {\cal MSV}(X)\to {\cal MSV}(X)$ 
is defined but is not associative. Rather, for every open set $U$ the space 
of sections of ${\cal MSV}$ over $U$ forms a vertex algebra and sections 
of $\cal O$ over $U$ are mapped to the set of pairwise commuting bosonic 
fields in $\Gamma(U,{\cal MSV})$. This type of sheaf was called 
a quasi-loop-coherent sheaf of vertex algebras in \cite{borvert}. 

One then observes that cohomology ${\rm H^*}(X,{\cal MSV})$
of the chiral de Rham complex is provided with a natural structure of 
vertex algebra, essentially via a cup product, see \cite{borvert} for
details. This vertex algebra has a natural N=2 structure, if $X$ is
Calabi-Yau. The corresponding structures of topological algebras, see 
Section 2, should correspond to the A and B models. In fact, it was 
the topological twist of the above algebra that was considered in
\cite{MSV}, so their definition of $L$ is slightly different.
As a result, the (half-integer) notations for the fermionic modes 
of the operators $\varphi$ and $\psi$ that we have used above are 
different from the (integer) notations used in \cite{MSV}. 
However, in terms of the natural modes $\varphi_{(n)}$ and $\psi_{(n)}$ 
our notations coincide.

It is worth mentioning that ${\cal MSV}(X)$ possesses a natural filtration,
such that the graded object is a quasi-coherent sheaf isomorphic
to a tensor product of an infinite number of copies of symmetric and exterior
algebras of the tangent and cotangent sheaves on $X$. This remark
allows one to show that the elliptic genus of the variety $X$ can
be naturally formulated in terms of the supertrace over the cohomology
of ${\cal MSV}(X)$ of the operator $y^{J[0]}q^{L[0]}$. The reader is referred
to \cite{borlibg} for details. 

It would be very interesting to compare the approach of \cite{MSV} with
that of the monograph of Tamanoi \cite{tamanoi}.

\section{Vertex algebras of Calabi-Yau hypersurfaces in toric varieties}
We will now talk about the contents of the paper \cite{borvert}.
Its major achievement is an explicit calculation of the cohomology 
of the chiral de Rham complex for a generic Calabi-Yau hypersurface $X$ 
in a smooth toric nef-Fano variety $\PP$. There is also some progress
made in the problem of defining chiral de Rham complex for varieties
with Gorenstein singularities. The description uses certain vertex algebra
constructed from a lattice, whose definition will be provided below.

Before we can describe the cohomology of the chiral de Rham complex
of a generic Calabi-Yau toric hypersurface, we must recall the combinatorial 
data which define it that were discovered by Batyrev in \cite{bat.dual}.
Let $M_1$ and $N_1$ be two dual lattices of rank $d+1$ that contain dual 
reflexive polyhedra $\Delta_1$ and $\Delta_1^*$. One defines dual lattices 
$M=M_1\oplus \ZZ$ and $N=N_1\oplus \ZZ$ of rank $d+2$ and cones 
$K=\{(t\Delta_1,t),\,t\geq 0\}$ and $K^*=\{(t\Delta_1^*,t),\,t\geq 0\}$ in
$M$ and $N$ respectively. The conditions on $\Delta_1$ and $\Delta_1^*$
to form a dual pair of reflexive polytopes is that all their vertices
are lattice points, and that $K$ and $K^*$ are dual to each other, as 
the notation suggests. To specify a nef-Fano toric variety $\PP$ 
one also chooses a fan $\Sigma$ in $N_1$ which subdivides the minimum
fan defined by the faces of $\Delta_1^*$. Toric variety $\PP$ is smooth
if and only if all cones of $\Sigma$ are generated by a part of the 
basis of the lattice $N_1$. 

We denote the bilinear form on $M\oplus N$ by $\cdot$. We also denote
$(0,1)\in M$ by $\deg$ and $(0,1)\in N$ by $\deg^*$. We call $\deg\cdot n$ 
and $\deg^*\cdot m$ the degree of $n\in N$ and $m\in M$ respectively.
Codimension one polytopes $\Delta_1+\deg$ and $\Delta_1^*+\deg^*$
are denoted by $\Delta$ and $\Delta^*$ respectively. A generic hypersurface 
$X$ in $\PP$ is defined by a generic collection of coefficients $f_m$ for 
all lattice points $m\in \Delta$. It is (in general only partial) 
desingularization of ${\rm Proj}(\CC[K]/f)$. Here $f$ is an element of degree
one defined by $\sum_{m\in\Delta} f_mx^m$ where $x$ is a dummy variable
used to write $\CC[K]$ in a multiplicative form.

To explain the results of \cite{borvert}, we need to introduce
a vertex algebra associated to the lattice $M \oplus N$. As a vector
space, it is isomorphic to the vector space of the vertex algebra
of $2(d+2)$ free fermions and $2(d+2)$ free bosons constructed in 
Section 4 tensored with $\CC[M\oplus N]$. For any subset $I$ of 
$M\oplus N$ we denote by $\Fock_I$ the space obtained by tensoring 
(over $\CC$) of the vertex algebra of Section 4 and vector
space $\CC[I]$. 

First of all, we need to define how the fields of 
$\Fock_{0\oplus 0}$ act on $\Fock_{M\oplus N}$. 
We use the notations of \cite{borvert} and denote the bosonic fields 
of $\Fock_{0\oplus 0}$ by $m\cdot B (z)$ and $n\cdot A(z)$ and fermionic 
fields by $m\cdot \Phi(z)$ and $n\cdot \Psi(z)$. Here $A,B,\Phi,\Psi$ are 
vector 
valued, and $m$ and $n$ are in $M\otimes\CC$ and $n\otimes \CC$ respectively.
For a fixed pair of lattice elements $(m,n)$, the action of fermionic
fields $\Phi$ and $\Psi$  on $Fock_{m\oplus n}$ is simply induced 
by their action on $\Fock_{0\oplus 0}$. The action of the bosonic fields
is modified so that the zero modes $(m_1\cdot B)[0]$ and $(n_1\cdot A)[0]$
do not annihilate $\Fock_{m\oplus n}$. Rather, they act by a scalar 
multiplication by $m_1\cdot n$ and $n_1\cdot m$ respectively.

To define state-field correspondence we still need to specify which 
fields of $\Fock_{M\oplus N}$ correspond to the elements of $\Fock_{m\oplus n}$
with non-zero $(m,n)$.  We denote the element 
$(\vac,m\oplus n)\in \Fock_{M\oplus N}$ by $|m,n\rangle$, and our
first goal is to describe the field $Y(|m,n\rangle,z)$.
For all $(m_1,n_1)\in M\oplus N$ all modes of the field 
$Y(|m,n\rangle,z)$ map $\Fock_{m_1\oplus n_1}$ to $\Fock_{(m+m_1)\oplus
(n+n_1)}$. We denote the endomorphism of $\Fock_{M\oplus N}$ that 
commutes with all non-zero modes of $A$ and $B$ and sends $|m_1,n_1\rangle$
to $|m+m_1,n+n_1\rangle$ by $\gamma_{m,n}$.
Then 
$$
Y(|m,n\rangle,z):=  \gamma_{m,n}
c_{m,n}\no\ee^{\sum_{0\neq k\in \ZZ} (-z^{-k}/k)((m\cdot B)[k]+(k\cdot A)[k])}
\on z^{m\cdot B[0]} z^{n\cdot A[0]}
$$
where $c_{m,n}$ acts on $\Fock_{m_1,n_1}$ by multiplication by 
$(-1)^{m\cdot n_1}$.
One defines $Y(a,z)$ for other elements of $\Fock_{m\oplus n}$ by inserting
appropriate free fields and their derivatives inside the normal ordering.
It is not hard to see that these operators are well defined, moreover
one can show that they are mutually local and satisfy some nice OPEs.
Instead of the complicated notation above we use 
$$
Y(|m,n\rangle,z):= \no\ee^{\int (m\cdot B+n\cdot A)(z)\,dz} \on
$$
and similarly for other elements of $\Fock_{m\oplus n}$.
All $|m,n\rangle$ are bosonic and satisfy 
$$
Y(|m,n\rangle,z)Y(|m_1,n_1\rangle,w)=\frac 
{\no \ee^{\int (m\cdot B+n\cdot A)(z)\,dz} 
\ee^{\int (m_1\cdot B+n_1\cdot A)(w)\,dw} \on}
{(z-w)^{m\cdot n_1+m_1\cdot n}}
$$
which can be Taylor expanded around $z=w$ to give the OPEs.
The details of this calculation could be found in \cite{Kac}.
In general, it is straightforward to calculate OPEs in the lattice
algebra, but the resulting expressions could be quite complicated.

The vertex algebra $\Fock_{M\oplus N}$ can be equipped with the following 
N=2 structure of central charge $\hat c= d$. Notice that the rank of 
the lattice $M$ is $d+2$. We will call it Calabi-Yau N=2 structure,
because we will see in a second that it is related to the N=2 structure
of the cohomology of the chiral de Rham complex of Calabi-Yau hypersurfaces
in toric varieties.
$$
G^+_{CY}(z):=(A\cdot\Phi)(z)-{\rm deg}\cdot\partial_z \Phi(z)
$$
$$
G^-_{CY}(z):= (B\cdot\Psi)(z)-{\rm deg}^*\cdot\partial_z \Psi(z)
$$
$$
J_{CY}(z):=\,\no(\Phi\cdot\Psi)(z)\on +\, {\rm deg}\cdot B(z)
-{\rm deg}^*\cdot A(z)
$$
$$
L_{CY}(z):=\,\no (B\cdot A)(z)\on +
\frac 12\no (\partial_z\Phi\cdot\Psi-\Phi\cdot \partial_z\Psi)(z)\on
-\frac 12{\rm deg}^*\cdot\partial_z A(z)
-\frac 12{\rm deg}\cdot\partial_z B(z)
.
$$

We will also need to describe a certain deformation of the vertex algebra
structure on $\Fock_{M\oplus N}$ defined by the fan $\Sigma_1$ 
which is used to define the ambient toric variety $\PP$. Fan $\Sigma_1$
naturally gives rise to a (generalized) fan $\Sigma$ in $N$ by simply 
extending all of the cones of $\Sigma_1$ in the direction of $\deg^*$. 
One can then define a vertex algebra structure on 
$\Fock_{M\oplus N}$ by changing $\gamma_{m,n}$ as follows. 
$$\gamma_{m,n}|m_1,n_1\rangle = 
\left\{
\begin{array}{ll}
|m+m_1,n+n_1\rangle & {\rm if~there~exists~}C\in\Sigma,~{\rm such~that~}
n_1,n\in C\\
0 &{\rm otherwise}
\end{array}
\right.
$$
This gives a different structure of the vertex algebra, and we denote it
by $\Fock_{M\oplus N}^{\Sigma}$ to distinguish it from the usual
structure on $\Fock_{M\oplus N}$. Operator product expansions of 
the new fields are either identical to the OPEs in $\Fock_{M\oplus N}$
or vanish. The N=2 superconformal structure on this new algebra is 
given by the same formulas. 

The following theorem is the main result of \cite{borvert} in the 
case of smooth ambient toric variety $\PP$.
\begin{theo}
{\rm 
Let $\Delta$ and $\Delta^*$ be as above, and let $\Sigma_1$ define
a non-singular toric variety $\PP$.  Let $f:(\Delta\cap M)\to \CC$ 
define a generic hypersurface $X\subset \PP$. We pick a generic 
collection of coefficients  $\{g_n,n\in\Delta^*\cap N\}$. Then 
the cohomology of the chiral de Rham complex of $X$ is isomorphic
as a vertex algebra to the BRST quotient of $\Fock_{M\oplus N}^\Sigma$ 
by the BRST operator
$${\cal BRST}_{f,g}:=\oint(
\sum_{m\in \Delta} f_m
(m\cdot \Phi)(z) \no\ee^{\int
m\cdot B(z)}\on+
\sum_{n\in \Delta^*} g_n
(n\cdot \Psi)(z) \no\ee^{\int
n\cdot A(z)}\on
)dz.$$ 
Moreover, the N=2 superconformal structure on ${\rm H}^*({\cal MSV}(X))$
coincides with the structure induced by $G^\pm_{CY},J_{CY},L_{CY}$
introduced above.
}
\label{mainborvert}
\end{theo}

We remark that BRST operator above supercommutes with the fields of
N=2 algebra which allows one to induce these fields on the BRST
quotient. Moreover, it was shown in \cite{borvert} that eigenvalues of $L[0]$
are non-negative on the BRST cohomology, even though they can definitely
be negative on elements of $\Fock_{M\oplus N}^\Sigma$. In addition, all
eigenspaces of $L[0]$ have finite dimension.

Notice that a nice feature of the above result is that the cohomology 
of the chiral de Rham complexes for mirror symmetric toric hypersurfaces
are obviously related to each other as deformations of a single family 
of the vertex algebras where one uses $\Fock_{M\oplus N}$ instead of
$\Fock_{M\oplus N}^\Sigma$. Then the only difference between mirror
pictures is that one switches the roles of $M$ and $N$ which amounts
exactly to the mirror involution $G^\pm\to G^\mp,J\to-J,L\to L$.
Besides, when $\PP=\PP^4$ and one picks $f_m$ and $g_n$ to be non-zero
at the vertices of the simplices only, one appears to recover the 
(finite quotient of) tensor product of vertex algebras with fractional 
charges that corresponds to the Gepner model \cite{gepner}.
Some details of this correspondence are still to be worked out, but
note the paper \cite{feigin.semikhatov}.

One also observes that combinatorial description of the N=2 superconformal
vertex algebra as a BRST quotient of a certain lattice vertex algebra 
makes perfect sense whether or not the ambient toric variety $\PP$ is 
non-singular. This prompts one to try to define chiral de Rham complex for 
a singular Calabi-Yau hypersurface $X={\rm Proj}(\CC[K]/f)$ according
to the general philosophy that one should be able to understand mirror
symmetry without using partial crepant desingularizations, see for
example \cite{bat.bor}. On the other hand it suggests that this chiral
de Rham complex will depend not only on the scheme structure of $X$ 
but also on some mysterious coefficients $g_n$. In the smooth case these
coefficients are rather irrelevant (all one needs is for them to be non-zero)
but in general they seem to be of importance. 

Paper \cite{borvert} contains a definition of a sheaf of vertex algebras 
${\cal MSV}(X)$ for a Calabi-Yau hypersurface in a toric variety. Namely,
for a toric affine chart $\AA_C$ of $\PP$ that corresponds to a cone 
$C_1\in\Sigma_1$ one considers a subcone $C\subset K$ defined 
by $\{(n_1,t)\in K,~{\rm s.t.~} n_1\in C\}$. Then sections of 
${\cal MSV}(X)$ over the intersection of $X$ with $\AA_C$ are 
defined as BRST quotient of $\Fock_{M\oplus C}$ by ${\rm BRST}_{f,g}$.
Here one induces the vertex algebra structure from $\Fock_{M\oplus N}$
to $\Fock_{M\oplus C}$ and ignores $n\notin C$ in the definition 
of ${\rm BRST}_{f,g}$. Some progress is made in the singular case,
but the exact analog of \ref{mainborvert} is still an open problem.
However, we were able to use the results of this analysis in the singular
case to define the elliptic genus of a singular Calabi-Yau hypersurface 
in a toric variety and prove that it satisfies the expected mirror duality,
see \cite{borlibg}.

\section{Open questions}
In this short section we briefly describe major open problems 
and minor technical obstacles that are still to be faced in the vertex algebra
approach to mirror symmetry.

The most important problem is that it is not clear how to see instanton
corrections in terms of the cohomology of the chiral de Rham complex
of a variety $X$. Geometrically, chiral de Rham complex seems to deal 
with the neighborhood of the constant loops in the loop space of $X$, 
while the instanton corrections are more global in nature. One may have 
to introduce some modules over the vertex algebra 
${\rm H}^*({\cal MSV}(X))$ to deal with this difficulty. This is
also related to the problem of putting together holomorphic and 
antiholomorpic parts of the N=(2,2) superconformal field theory associated
to a Calabi-Yau manifold. One distinct possibility, which the author
plans to explore, is that the true vertex algebra of the superconformal theory
associated to a generic Calabi-Yau hypersurface in a toric variety 
is BRST quotient of $\Fock_{M\oplus N}$, and that $\Fock_{M\oplus N}^\Sigma$
appears when one expands the correlators around the limiting point 
that corresponds to the degeneration of $\CC[M\oplus N]$ into
$\CC[M\oplus N]^\Sigma$. Then one hopes to recover instanton corrections
as higher terms in this expansion.

Another big issue is a possible extension of these definitions to 
curves of higher genus. This is always a highly non-trivial problem
in conformal field theory, and it is interesting to see if an explicit
description could be obtained in the case of hypersurfaces in toric 
varieties.

Less formidable problems include the extension of all definitions 
and results to the case of singular varieties with some mild singularities.
One also wants to show that the families of vertex algebras given
by Theorem \ref{mainborvert} are flat in the appropriate sense, that 
is the dimensions of the $L[0]$ graded components are generically constant.
One should somehow see how GKZ hypergeometric system of differential 
equations appears in the context of the above vertex algebra. Also, it
was suggested by Yi-Zhi Huang that one should expect this vertex algebra 
to have an invariant bilinear form in the sense of \cite{FHL}.

It is our hope that the interplay of vertex algebras and 
algebraic geometry will enrich both fields and will provide deep mathematical
understanding of conformal field theories that are currently only defined in
terms of path integrals.

\end{document}